\documentclass[a4paper]{amsart}

\usepackage[mathscr]{eucal}
\usepackage{amsmath,amsfonts,amssymb}
\usepackage{amsthm}
\usepackage{mathrsfs}
\usepackage{bm}
\usepackage{cite}
\usepackage[dvips]{graphicx,color}
\usepackage{fancyhdr}

\newcommand{\mathsym}[1]{{}}

\theoremstyle{plain}

\newtheorem{theorem}{Theorem}[section]

\newtheorem{lemma}[theorem]{Lemma}

\theoremstyle{definition}
\newtheorem{definition}{Definition}[section]

\begin{document}
\title[On $3$-lattices and spherical designs]{On $3$-lattices and spherical designs}
\author{Junichi Shigezumi}
\date{October 23, 2008}

\maketitle \vspace{-0.1in}
\begin{center}
Graduate School of Mathematics Kyushu University\\
Hakozaki 6-10-1 Higashi-ku, Fukuoka, 812-8581 Japan\\
{\it E-mail address} : j.shigezumi@math.kyushu-u.ac.jp \vspace{-0.05in}
\end{center} \quad

\begin{quote}
{\small\bfseries Abstract.}
An integral lattice which is generated by some vectors of norm $q$ is called $q$-lattice. Classification of $3$-lattices of dimension at most four is given by Mimura (On $3$-lattice, 2006). As a expansion, we give a classification of $3$-lattices of dimension at most seven. In addition, we consider the spherical designs from its shells.\\  \vspace{-0.15in}

\vspace{-0.15in}

 \vspace{0.15in}
\end{quote}

\section{Introduction}

Let $L$ be an Euclidean lattice. We call the squared norm of a vector of the lattice {\it norm of the vector}. Then, the set $s_m (L)$ of vectors of the lattice $L$ which has the same norm $m$ is called the {\it shell of the lattice}, i.e. $s_m (L) := \{ x \in L \: ; \: (x, x) = m \}$. Moreover, we call the shell of {\it minimum} $\min_{x \in L \setminus \{ 0 \}} (x, x)$ of the lattice $L$ the {\it minimal shell}. If an integral lattice $L$ is generated by the vectors of $s_q (L)$, then we call $L$ a $q$-lattice.

Let $S^{d-1} := \{ (x_1, \ldots, x_d) \in \mathbb{R}^d \, ; \, x_1^2 + \cdots + x_d^2 = 1 \}$ be the Euclidean sphere for $d \geqslant 1$.

\begin{definition}[Spherical design \cite{DGS}]\label{def-design}
Let $X$ be a non-empty finite set on the Euclidean sphere $S^{d-1}$, and let $t$ be a positive integer. $X$ is called a spherical $t$-design if
\begin{equation}
\frac{1}{| S^{d-1} |} \int_{S^{d-1}} f(\xi) \, d \xi \ = \ \frac{1}{| X |} \hspace{0.05in} \sum_{\xi \in X} \hspace{0.05in} f(\xi) \label{eq-design}
\end{equation}
for every polynomial $f(x) = f(x_1, \ldots, x_d)$ of degree at most $t$.
\end{definition}

Here, the left-hand side of the above equation (\ref{eq-design}) means the average on the sphere $S^{d-1}$, and the right-hand side means the average on the finite subset $X$. Thus, if $X$ is a spherical design, then $X$ gives a certain approximation of the sphere $S^{d-1}$.

For every nonempty shell $s_m (L)$ of lattice $L$, we consider normalization $\frac{1}{\sqrt{m}} s_m (L)$, which is a finite set on a Euclidean sphere. Then, we can consider spherical designs from shells of lattices. A lattice whose minimal shell is a spherical $5$-design is said to be {\it strongly perfect}.

Now, we have the following facts (We omit details on definitions. Please refer to the references.):

\begin{theorem}[see \cite{V} or \cite{P}]
A strongly perfect integral lattice of minimum $2$ is isometric to one of root lattices $\mathbb{A}_1$, $\mathbb{A}_2$, $\mathbb{D}_4$, $\mathbb{E}_6$, $\mathbb{E}_7$, or $\mathbb{E}_8$. Furthermore, the minimal shell of the lattice is an $7$-design only for the case of the lattice $\mathbb{E}_8$.\label{th-m2-t5}
\end{theorem}

\begin{theorem}[Venkov \cite{V}, Theorem 7.4]
The strongly perfect lattices that are integral and of minimum $3$ are $O_1$, $O_7$, $O_{16}$, $O_{22}$, and $O_{23}$. Furthermore, the minimal shell is an $7$-design only for the case of the lattice $O_{23}$.\label{th-m3-t5}
\end{theorem}

Note that all the lattices in Theorem \ref{th-m2-t5} (resp. Theorem \ref{th-m3-t5}) are $2$-lattices (resp. $3$-lattices). To classify $q$-lattices seems to be one of the ways to look for `good' spherical designs.

Now, we have the following theorem proved by Mimura:

\begin{theorem}[Mimura \cite{M}]
Let $L$ be a Euclidean lattice, which has orthogonal decomposition into irreducible lattices
\begin{equation*}
L = \mathbb{Z}^m \perp L_1 \perp \cdots \perp L_n.
\end{equation*}
Then, $L$ is a $3$-lattice if and only if every $L_i$ ($1 \leqslant i \leqslant n$) is $2$-lattice, $3$-lattice of minimum $2$, or $3$-lattice of minimum $3$, and satisfies one of the following conditions:
\begin{enumerate}
\item $m = 0$, and every $L_i$ is $3$-lattice of minimum $2$ or minimum $3$.

\item $m = 1$ or $2$, and one of $L_i$ is $2$-lattice or $3$-lattice of minimum $2$. Furthermore, if one of $L_i$ is a $1$ dimensional $2$-lattice, then one of $L_i$ is at least $2$ dimensional $2$-lattice or $3$-lattice of minimum $2$.

\item $m = 3$, and one of $L_i$ is $2$-lattice or $3$-lattice of minimum $2$.

\item $m \geqslant 4$.
\end{enumerate}
\end{theorem}

In addition, irreducible $2$-lattices are isometric to one of root lattices $\mathbb{A}_n$ for $n \geqslant 1$, $\mathbb{D}_n$ for $n \geqslant 4$, or $\mathbb{E}_n$ for $n = 6, 7, 8$. Thus, to classify $3$-lattices, we have only to classify irreducible $3$-lattices of minimum $2$ or minimum $3$.

We classified irreducible $3$-lattices of minimum $2$ or minimum $3$ and of at most $7$ dimension by numerical calculation with `Magma'. We saw the result at the end of this note, in the appendix. The following table is the numbers of $d$ dimensional irreducible $3$-lattices of minimum $2$ and minimum $3$. (cf. \cite{S})

\begin{table}[h]
\begin{tabular}{cccccccc}
\hline
 & $1$ & $2$ & $3$ & $4$ & $5$ & $6$ & $7$\\
\hline
min. $2$ & $-$ & $1$ & $4$ & $18$ & $90$ & $591$ & $5449$\\
min. $3$ & $1$ & $1$ & $3$ & $9$ & $39$ & $247$ & $2446$\\
\hline\\
\end{tabular}
\caption{the numbers of irreducible $3$-lattices}
\end{table}

After that, we calculate spherical designs of $3$-norm shells of all the $3$-lattices, which are reducible and irreducible. In the following section, we pick up all the $3$-lattices whose shells of norm $3$ are at least $3$-designs.

\setcounter{tocdepth}{1}
\tableofcontents

\section{Spherical designs from $3$-lattices}

Let $X$ be a nonempty finite set on the Euclidean sphere $S^{d-1}$ ($\subset \mathbb{R}^{d}$). We denote the distance set of $X$ by $A(X) := \{ (x, y) ; x, y \in X, x \ne y \}$, then we call $X$ a $s$-distance set if $| A(X) | = s$. Now, $X$ is said to be a $(d, n, s, t)$-configuration if $X \subset S^{d-1}$ is of order $n (:= | X |)$, a $s$-distance set, and a spherical $t$-design.

We denote by $GM(L)$ a Gram matrix of a lattice $L$, and by $\det (L) := \det (GM(L))$ the determinant of a lattice. Furthermore, we denote $(L)^n := L \oplus \cdots \oplus L$.\\

\subsection{$d = 1$} We have only one $3$-lattice $\sqrt{3} \mathbb{Z}$. Since we have $S^0 = \frac{1}{\sqrt{m}} s_m (\sqrt{3} \mathbb{Z})$, every shell of this lattice is a spherical $t$-design for all the positive integer $t$. Then, for every $1$-dimensional lattice, we describe $(d, n, s, t) = (1, 2, 1, *)$.

\quad

\subsection{$d = 2$} We have only one $3$-lattice whose shell of norm $3$ is $3$-design, which is $(\sqrt{3} \mathbb{Z})^2$. Here, we have $\det ((\sqrt{3} \mathbb{Z})^2) = 9$ and $(d, n, s, t) = (2,4,2,3)$.

Now, the following table is $(d, n, s, t)$-configuration of each shell of norm $m \leqslant 12$ of the lattice:
\begin{center}
\begin{tabular}{r}
\hline
$m$\\
\hline
$3$\\
$6$\\
$12$\\
\hline
\end{tabular}
\begin{tabular}{cccc}
\hline
$d$ & $n$ & $s$ & $t$\\
\hline
$2$ & $4$ & $2$ & $3$\\
$2$ & $4$ & $2$ & $3$\\
$2$ & $4$ & $2$ & $3$\\
\hline
\end{tabular}
\end{center}

\newpage

\subsection{$d = 3$} We have three $3$-lattices whose shells of norm $3$ are $3$-design, which are $L_{3,3}$, $(\sqrt{3} \mathbb{Z})^3$, and $L_{3,1}$. A Gram matrix of each lattice is the following:

\begin{equation*}
GM(L_{3,3}) = \left[
\begin{matrix}
 3 & 1 & 1 \\
 1 & 3 & -1 \\
 1 & -1 & 3
\end{matrix}
\right],
 \quad
GM(L_{3,1}) = \left[
\begin{matrix}
 1 & 0 & 0 \\
 0 & 2 & 1 \\
 0 & 1 & 2
\end{matrix}
\right].
\end{equation*}
Here, we have $\det (L_{3,3}) = 16$, $\det ((\sqrt{3} \mathbb{Z})^3) = 27$, and $\det (L_{3,1}) = 3$.

Now, the following table is $(d, n, s, t)$-configuration of each shell of norm $m \leqslant 12$ of each lattice:
\begin{center}
\begin{tabular}{r}
\quad\\
\hline
$m$\\
\hline
$1$\\
$2$\\
$3$\\
$4$\\
$5$\\
$6$\\
$7$\\
$8$\\
$9$\\
$10$\\
$11$\\
$12$\\
\hline
\end{tabular}
\quad
\begin{tabular}{cccc}
\multicolumn{4}{c}{$L_{3,3}$}\\
\hline
$d$ & $n$ & $s$ & $t$\\
\hline
\ & \ & \ & \\
\ & \ & \ & \\
$3$ & $8$ & $3$ & $3$\\
$3$ & $6$ & $2$ & $3$\\
\ & \ & \ & \\
\ & \ & \ & \\
\ & \ & \ & \\
$3$ & $12$ & $4$ & $3$\\
\ & \ & \ & \\
\ & \ & \ & \\
$3$ & $24$ & $9$ & $3$\\
$3$ & $8$ & $3$ & $3$\\
\hline
\end{tabular}
\quad
\begin{tabular}{cccc}
\multicolumn{4}{c}{$(\sqrt{3} \mathbb{Z})^3$}\\
\hline
$d$ & $n$ & $s$ & $t$\\
\hline
\ & \ & \ & \\
\ & \ & \ & \\
$3$ & $6$ & $2$ & $3$\\
\ & \ & \ & \\
\ & \ & \ & \\
$3$ & $12$ & $4$ & $3$\\
\ & \ & \ & \\
\ & \ & \ & \\
$3$ & $8$ & $3$ & $3$\\
\ & \ & \ & \\
\ & \ & \ & \\
$3$ & $6$ & $2$ & $3$\\
\hline
\end{tabular}
\quad
\begin{tabular}{cccc}
\multicolumn{4}{c}{$L_{3,1}$}\\
\hline
$d$ & $n$ & $s$ & $t$\\
\hline
$1$ & $2$ & $1$ & $*$\\
$2$ & $6$ & $3$ & $1$\\
$3$ & $12$ & $6$ & $3$\\
$1$ & $2$ & $1$ & $*$\\
\ & \ & \ & \\
$3$ & $18$ & $8$ & $1$\\
$3$ & $12$ & $7$ & $1$\\
$2$ & $6$ & $3$ & $1$\\
$3$ & $14$ & $7$ & $1$\\
$3$ & $12$ & $7$ & $1$\\
$3$ & $12$ & $7$ & $1$\\
$3$ & $12$ & $6$ & $3$\\
\hline
\end{tabular}
\end{center}

\quad\\

\subsection{$d = 4$} We have two $3$-lattices whose shells of norm $3$ are $3$-designs, which are $(\sqrt{3} \mathbb{Z})^4$ and $\mathbb{Z}^4$. Here, we have $\det ((\sqrt{3} \mathbb{Z})^4) = 81$ and $\det (\mathbb{Z}^4) = 1$.

Now, the following table is $(d, n, s, t)$-configuration of each shell of norm $m \leqslant 12$ of each lattice:
\begin{center}
\begin{tabular}{r}
\quad\\
\hline
$m$\\
\hline
$1$\\
$2$\\
$3$\\
$4$\\
$5$\\
$6$\\
$7$\\
$8$\\
$9$\\
$10$\\
$11$\\
$12$\\
\hline
\end{tabular}
\quad
\begin{tabular}{cccc}
\multicolumn{4}{c}{$(\sqrt{3} \mathbb{Z})^4$}\\
\hline
$d$ & $n$ & $s$ & $t$\\
\hline
\ & \ & \ & \\
\ & \ & \ & \\
$4$ & $8$ & $2$ & $3$\\
\ & \ & \ & \\
\ & \ & \ & \\
$4$ & $24$ & $4$ & $5$\\
\ & \ & \ & \\
\ & \ & \ & \\
$4$ & $32$ & $6$ & $3$\\
\ & \ & \ & \\
\ & \ & \ & \\
$4$ & $24$ & $4$ & $5$\\
\hline
\end{tabular}
\quad
\begin{tabular}{cccc}
\multicolumn{4}{c}{$\mathbb{Z}^4$}\\
\hline
$d$ & $n$ & $s$ & $t$\\
\hline
$4$ & $8$ & $2$ & $3$\\
$4$ & $24$ & $4$ & $5$\\
$4$ & $32$ & $6$ & $3$\\
$4$ & $24$ & $4$ & $5$\\
$4$ & $48$ & $10$ & $3$\\
$4$ & $96$ & $12$ & $5$\\
$4$ & $64$ & $14$ & $3$\\
$4$ & $24$ & $4$ & $5$\\
$4$ & $104$ & $18$ & $3$\\
$4$ & $144$ & $20$ & $5$\\
$4$ & $96$ & $20$ & $3$\\
$4$ & $96$ & $12$ & $5$\\
\hline
\end{tabular}
\end{center}

\newpage

\subsection{$d = 5$} We have three $3$-lattices whose shells of norm $3$ are $3$-design, which are $L_{5,3}$, $(\sqrt{3} \mathbb{Z})^5$, and $\mathbb{Z}^5$. A Gram matrix of $L_{5,3}$ is
\begin{equation*}
GM(L_{5,3}) = \left[
\begin{matrix}
 3 & 1 & 1 & 1 & 1 \\
 1 & 3 & 1 & 1 & 1 \\
 1 & 1 & 3 & 1 & -1 \\
 1 & 1 & 1 & 3 & -1 \\
 1 & 1 & -1 & -1 & 3
\end{matrix}
\right].
\end{equation*}
Here, we have $\det (L_{5,3}) = 48$, $\det ((\sqrt{3} \mathbb{Z})^5) = 243$, and $\det (\mathbb{Z}^5) = 1$.

Now, the following table is $(d, n, s, t)$-configuration of each shell of norm $m \leqslant 12$ of each lattice:
\begin{center}
\begin{tabular}{r}
\quad\\
\hline
$m$\\
\hline
$1$\\
$2$\\
$3$\\
$4$\\
$5$\\
$6$\\
$7$\\
$8$\\
$9$\\
$10$\\
$11$\\
$12$\\
\hline
\end{tabular}
\quad
\begin{tabular}{cccc}
\multicolumn{4}{c}{$L_{5,3}$}\\
\hline
$d$ & $n$ & $s$ & $t$\\
\hline
\ & \ & \ & \\
\ & \ & \ & \\
$5$ & $20$ & $3$ & $3$\\
$5$ & $30$ & $4$ & $3$\\
\ & \ & \ & \\
\ & \ & \ & \\
$5$ & $60$ & $7$ & $3$\\
$5$ & $90$ & $8$ & $3$\\
\ & \ & \ & \\
\ & \ & \ & \\
$5$ & $180$ & $11$ & $3$\\
$5$ & $140$ & $12$ & $3$\\
\hline
\end{tabular}
\quad
\begin{tabular}{cccc}
\multicolumn{4}{c}{$(\sqrt{3} \mathbb{Z})^5$}\\
\hline
$d$ & $n$ & $s$ & $t$\\
\hline
\ & \ & \ & \\
\ & \ & \ & \\
$5$ & $10$ & $2$ & $3$\\
\ & \ & \ & \\
\ & \ & \ & \\
$5$ & $40$ & $4$ & $3$\\
\ & \ & \ & \\
\ & \ & \ & \\
$5$ & $80$ & $6$ & $3$\\
\ & \ & \ & \\
\ & \ & \ & \\
$5$ & $90$ & $8$ & $3$\\
\hline
\end{tabular}
\quad
\begin{tabular}{cccc}
\multicolumn{4}{c}{$\mathbb{Z}^5$}\\
\hline
$d$ & $n$ & $s$ & $t$\\
\hline
$5$ & $10$ & $2$ & $3$\\
$5$ & $40$ & $4$ & $3$\\
$5$ & $80$ & $6$ & $3$\\
$5$ & $90$ & $8$ & $3$\\
$5$ & $112$ & $10$ & $3$\\
$5$ & $224$ & $12$ & $3$\\
$5$ & $320$ & $14$ & $3$\\
$5$ & $200$ & $16$ & $3$\\
$5$ & $250$ & $18$ & $3$\\
$5$ & $560$ & $20$ & $3$\\
$5$ & $560$ & $22$ & $3$\\
$5$ & $400$ & $24$ & $3$\\
\hline
\end{tabular}
\end{center}

\newpage

\subsection{$d = 6$}\quad

\subsubsection{$3$-lattices of minimum $3$} We have four $3$-lattices whose shells of norm $3$ are $3$-design, which are $L_{6,3,1}$, $L_{6,3,2}$, $(L_{3,3})^2$, and $(\sqrt{3} \mathbb{Z})^6$. A Gram matrix of each lattice is
\begin{equation*}
GM(L_{6,3,1}) = \left[
\begin{matrix}
 3 & 1 & 1 & 1 & 1 & 1 \\
 1 & 3 & 1 & 1 & 1 & 1 \\
 1 & 1 & 3 & 1 & 1 & 1 \\
 1 & 1 & 1 & 3 & 1 & -1 \\
 1 & 1 & 1 & 1 & 3 & -1 \\
 1 & 1 & 1 & -1 & -1 & 3
\end{matrix}
\right],
 \quad
GM(L_{6,3,2}) = \left[
\begin{matrix}
 3 & 1 & 1 & 1 & 1 & 0 \\
 1 & 3 & 1 & 0 & -1 & 1 \\
 1 & 1 & 3 & -1 & 0 & -1 \\
 1 & 0 & -1 & 3 & 1 & 1 \\
 1 & -1 & 0 & 1 & 3 & -1 \\
 0 & 1 & -1 & 1 & -1 & 3
\end{matrix}
\right].
\end{equation*}
Here, we have $\det (L_{6,3,1}) = 64$, $\det (L_{6,3,2}) = 125$, $\det ((L_{3,3})^2) = 256$, and $\det ((\sqrt{3} \mathbb{Z})^6) = 729$.

Now, the following table is $(d, n, s, t)$-configuration of each shell of norm $m \leqslant 12$ of each lattice:
\begin{center}
\begin{tabular}{r}
\quad\\
\hline
$m$\\
\hline
$1$\\
$2$\\
$3$\\
$4$\\
$5$\\
$6$\\
$7$\\
$8$\\
$9$\\
$10$\\
$11$\\
$12$\\
\hline
\end{tabular}
\quad
\begin{tabular}{cccc}
\multicolumn{4}{c}{$L_{6,3,1}$}\\
\hline
$d$ & $n$ & $s$ & $t$\\
\hline
\ & \ & \ & \\
\ & \ & \ & \\
$6$ & $32$ & $3$ & $3$\\
$6$ & $60$ & $4$ & $3$\\
\ & \ & \ & \\
\ & \ & \ & \\
$6$ & $192$ & $7$ & $3$\\
$6$ & $252$ & $8$ & $3$\\
\ & \ & \ & \\
\ & \ & \ & \\
$6$ & $480$ & $11$ & $3$\\
$6$ & $544$ & $12$ & $3$\\
\hline
\end{tabular}
\quad
\begin{tabular}{cccc}
\multicolumn{4}{c}{$L_{6,3,2}$}\\
\hline
$d$ & $n$ & $s$ & $t$\\
\hline
\ & \ & \ & \\
\ & \ & \ & \\
$6$ & $20$ & $4$ & $3$\\
$6$ & $30$ & $6$ & $3$\\
$6$ & $24$ & $6$ & $3$\\
$6$ & $60$ & $10$ & $3$\\
$6$ & $60$ & $10$ & $3$\\
$6$ & $60$ & $12$ & $3$\\
$6$ & $120$ & $16$ & $3$\\
$6$ & $144$ & $18$ & $3$\\
$6$ & $240$ & $20$ & $3$\\
$6$ & $200$ & $20$ & $3$\\
\hline
\end{tabular}
\quad
\begin{tabular}{cccc}
\multicolumn{4}{c}{$(L_{3,3})^2$}\\
\hline
$d$ & $n$ & $s$ & $t$\\
\hline
\ & \ & \ & \\
\ & \ & \ & \\
$6$ & $16$ & $4$ & $3$\\
$6$ & $12$ & $2$ & $3$\\
\ & \ & \ & \\
$6$ & $64$ & $6$ & $3$\\
$6$ & $96$ & $10$ & $3$\\
$6$ & $60$ & $4$ & $3$\\
\ & \ & \ & \\
\ & \ & \ & \\
$6$ & $240$ & $16$ & $3$\\
$6$ & $160$ & $6$ & $3$\\
\hline
\end{tabular}
\quad
\begin{tabular}{cccc}
\multicolumn{4}{c}{$(\sqrt{3} \mathbb{Z})^6$}\\
\hline
$d$ & $n$ & $s$ & $t$\\
\hline
\ & \ & \ & \\
\ & \ & \ & \\
$6$ & $12$ & $2$ & $3$\\
\ & \ & \ & \\
\ & \ & \ & \\
$6$ & $60$ & $4$ & $3$\\
\ & \ & \ & \\
\ & \ & \ & \\
$6$ & $160$ & $6$ & $3$\\
\ & \ & \ & \\
\ & \ & \ & \\
$6$ & $252$ & $8$ & $3$\\
\hline
\end{tabular}
\end{center}

\quad\\

\subsubsection{$3$-lattices of minimum $2$} We have three $3$-lattices whose shells of norm $3$ are $3$-design, which are $L_{6,2,1}$, $L_{6,2,2}$, and $L_{6,2,3}$. A Gram matrix of each lattice is
\begin{equation*}
L_{6,2,1}: \; \left[
\begin{matrix}
 2 & 0 & 0 & 0 & 0 & 1 \\
 0 & 2 & 0 & 0 & 0 & 1 \\
 0 & 0 & 2 & 0 & 0 & 1 \\
 0 & 0 & 0 & 2 & 0 & 1 \\
 0 & 0 & 0 & 0 & 2 & 1 \\
 1 & 1 & 1 & 1 & 1 & 3
\end{matrix}
\right],
 \quad
L_{6,2,2}: \; \left[
\begin{matrix}
 2 & 0 & 1 & 1 & 1 & 1 \\
 0 & 2 & 1 & 1 & 1 & 1 \\
 1 & 1 & 3 & 1 & 1 & 1 \\
 1 & 1 & 1 & 3 & 1 & 1 \\
 1 & 1 & 1 & 1 & 3 & 1 \\
 1 & 1 & 1 & 1 & 1 & 3
\end{matrix}
\right],
 \quad
L_{6,2,3}: \; \left[
\begin{matrix}
 2 & 1 & 1 & 1 & 1 & 1 \\
 1 & 3 & 1 & 1 & 1 & 1 \\
 1 & 1 & 3 & 0 & 0 & 0 \\
 1 & 1 & 0 & 3 & 0 & 0 \\
 1 & 1 & 0 & 0 & 3 & 0 \\
 1 & 1 & 0 & 0 & 0 & 3
\end{matrix}
\right].
\end{equation*}
Here, we have $\det (L_{6,2,1}) = 16$, $\det (L_{6,2,2}) = 64$, and $\det (L_{6,2,3}) = 81$.

Now, the following table is $(d, n, s, t)$-configuration of each shell of norm $m \leqslant 12$ of each lattice:
\begin{center}
\begin{tabular}{r}
\quad\\
\hline
$m$\\
\hline
$1$\\
$2$\\
$3$\\
$4$\\
$5$\\
$6$\\
$7$\\
$8$\\
$9$\\
$10$\\
$11$\\
$12$\\
\hline
\end{tabular}
\quad
\begin{tabular}{cccc}
\multicolumn{4}{c}{$L_{6,2,1}$}\\
\hline
$d$ & $n$ & $s$ & $t$\\
\hline
\ & \ & \ & \\
$6$ & $12$ & $2$ & $3$\\
$6$ & $64$ & $6$ & $3$\\
$6$ & $60$ & $4$ & $3$\\
\ & \ & \ & \\
$6$ & $160$ & $6$ & $3$\\
$6$ & $384$ & $14$ & $3$\\
$6$ & $252$ & $8$ & $3$\\
\ & \ & \ & \\
$6$ & $312$ & $10$ & $3$\\
$6$ & $960$ & $22$ & $3$\\
$6$ & $544$ & $12$ & $3$\\
\hline
\end{tabular}
\quad
\begin{tabular}{cccc}
\multicolumn{4}{c}{$L_{6,2,2}$}\\
\hline
$d$ & $n$ & $s$ & $t$\\
\hline
\ & \ & \ & \\
$2$ & $4$ & $2$ & $1$\\
$6$ & $32$ & $6$ & $3$\\
$6$ & $28$ & $4$ & $1$\\
\ & \ & \ & \\
$6$ & $96$ & $6$ & $3$\\
$6$ & $192$ & $14$ & $3$\\
$6$ & $124$ & $8$ & $1$\\
\ & \ & \ & \\
$6$ & $104$ & $10$ & $1$\\
$6$ & $480$ & $22$ & $3$\\
$6$ & $288$ & $12$ & $3$\\
\hline
\end{tabular}
\quad
\begin{tabular}{cccc}
\multicolumn{4}{c}{$L_{6,2,3}$}\\
\hline
$d$ & $n$ & $s$ & $t$\\
\hline
\ & \ & \ & \\
$1$ & $2$ & $1$ & $*$\\
$6$ & $24$ & $6$ & $3$\\
$5$ & $30$ & $5$ & $1$\\
$6$ & $40$ & $7$ & $1$\\
$6$ & $90$ & $10$ & $1$\\
$6$ & $24$ & $7$ & $1$\\
$6$ & $62$ & $9$ & $1$\\
$6$ & $280$ & $18$ & $1$\\
$5$ & $132$ & $13$ & $1$\\
$6$ & $120$ & $15$ & $1$\\
$6$ & $414$ & $22$ & $1$\\
\hline
\end{tabular}
\end{center}

\newpage

\subsubsection{$3$-lattices of minimum $1$} We have six $3$-lattices whose shells of norm $3$ are $3$-design, which are $\mathbb{Z}^6$, $\mathbb{Z}^2 \oplus \mathbb{D}_4$, $\mathbb{Z}^4 \oplus (\sqrt{2} \mathbb{Z})^2$, $\mathbb{Z}^2 \oplus \mathbb{A}_4$, $\mathbb{Z}^2 \oplus (\mathbb{A}_2)^2$, and $L_{6,1}$. A Gram matrix of $L_{6,1}$ is
\begin{equation*}
GM(L_{6,1}) = \left[
\begin{matrix}
 1 & 0 & 0 & 0 & 0 & 0 \\
 0 & 2 & 1 & 1 & 1 & 1 \\
 0 & 1 & 2 & 1 & 1 & 1 \\
 0 & 1 & 1 & 2 & 0 & 0 \\
 0 & 1 & 1 & 0 & 3 & 1 \\
 0 & 1 & 1 & 0 & 1 & 3
\end{matrix}
\right].
\end{equation*}
Here, we have $\det (\mathbb{Z}^6) = 1$, $\det (\mathbb{Z}^2 \oplus \mathbb{D}_4) = 4$, $\det (\mathbb{Z}^4 \oplus (\sqrt{2} \mathbb{Z})^2) = 4$, $\det (\mathbb{Z}^2 \oplus \mathbb{A}_4) = 5$, $\det (\mathbb{Z}^2 \oplus (\mathbb{A}_2)^2) = 9$, and $\det (L_{6,1}) = 16$.

Now, the following table is $(d, n, s, t)$-configuration of each shell of norm $m \leqslant 12$ of each lattice:
\begin{center}
\begin{tabular}{r}
\quad\\
\hline
$m$\\
\hline
$1$\\
$2$\\
$3$\\
$4$\\
$5$\\
$6$\\
$7$\\
$8$\\
$9$\\
$10$\\
$11$\\
$12$\\
\hline
\end{tabular}
\quad
\begin{tabular}{cccc}
\multicolumn{4}{c}{$\mathbb{Z}^6$}\\
\hline
$d$ & $n$ & $s$ & $t$\\
\hline
$6$ & $12$ & $2$ & $3$\\
$6$ & $60$ & $4$ & $3$\\
$6$ & $160$ & $6$ & $3$\\
$6$ & $252$ & $8$ & $3$\\
$6$ & $312$ & $10$ & $3$\\
$6$ & $544$ & $12$ & $3$\\
$6$ & $960$ & $14$ & $3$\\
$6$ & $1020$ & $16$ & $3$\\
$6$ & $876$ & $18$ & $3$\\
$6$ & $1560$ & $20$ & $3$\\
$6$ & $2400$ & $22$ & $3$\\
$6$ & $2080$ & $24$ & $3$\\
\hline
\end{tabular}
\quad
\begin{tabular}{cccc}
\multicolumn{4}{c}{$\mathbb{Z}^2 \oplus \mathbb{D}_4$}\\
\hline
$d$ & $n$ & $s$ & $t$\\
\hline
$2$ & $4$ & $2$ & $1$\\
$6$ & $28$ & $4$ & $1$\\
$6$ & $96$ & $6$ & $3$\\
$6$ & $124$ & $8$ & $1$\\
$6$ & $104$ & $10$ & $1$\\
$6$ & $288$ & $12$ & $3$\\
$6$ & $576$ & $14$ & $3$\\
$6$ & $508$ & $16$ & $1$\\
$6$ & $292$ & $18$ & $1$\\
$6$ & $728$ & $20$ & $1$\\
$6$ & $1440$ & $22$ & $3$\\
$6$ & $1056$ & $24$ & $3$\\
\hline
\end{tabular}
\quad
\begin{tabular}{cccc}
\multicolumn{4}{c}{$\mathbb{Z}^4 \oplus (\sqrt{2} \mathbb{Z})^2$}\\
\hline
$d$ & $n$ & $s$ & $t$\\
\hline
$4$ & $8$ & $2$ & $1$\\
$6$ & $28$ & $4$ & $1$\\
$6$ & $64$ & $6$ & $3$\\
$6$ & $124$ & $8$ & $1$\\
$6$ & $208$ & $10$ & $1$\\
$6$ & $288$ & $12$ & $3$\\
$6$ & $384$ & $14$ & $3$\\
$6$ & $508$ & $16$ & $1$\\
$6$ & $584$ & $18$ & $1$\\
$6$ & $728$ & $20$ & $1$\\
$6$ & $960$ & $22$ & $3$\\
$6$ & $1056$ & $24$ & $3$\\
\hline
\end{tabular}
\end{center}

\quad\\

\begin{center}
\begin{tabular}{r}
\quad\\
\hline
$m$\\
\hline
$1$\\
$2$\\
$3$\\
$4$\\
$5$\\
$6$\\
$7$\\
$8$\\
$9$\\
$10$\\
$11$\\
$12$\\
\hline
\end{tabular}
\quad
\begin{tabular}{cccc}
\multicolumn{4}{c}{$\mathbb{Z}^2 \oplus \mathbb{A}_4$}\\
\hline
$d$ & $n$ & $s$ & $t$\\
\hline
$2$ & $4$ & $2$ & $1$\\
$6$ & $24$ & $4$ & $1$\\
$6$ & $80$ & $6$ & $3$\\
$6$ & $114$ & $8$ & $1$\\
$6$ & $128$ & $10$ & $1$\\
$6$ & $260$ & $12$ & $1$\\
$6$ & $400$ & $14$ & $1$\\
$6$ & $424$ & $16$ & $1$\\
$6$ & $484$ & $18$ & $1$\\
$6$ & $688$ & $20$ & $1$\\
$6$ & $1040$ & $22$ & $1$\\
$6$ & $1040$ & $24$ & $1$\\
\hline
\end{tabular}
\quad
\begin{tabular}{cccc}
\multicolumn{4}{c}{$\mathbb{Z}^2 \oplus (\mathbb{A}_2)^2$}\\
\hline
$d$ & $n$ & $s$ & $t$\\
\hline
$2$ & $4$ & $2$ & $1$\\
$6$ & $16$ & $4$ & $1$\\
$6$ & $48$ & $6$ & $3$\\
$6$ & $88$ & $8$ & $1$\\
$6$ & $152$ & $10$ & $1$\\
$6$ & $204$ & $12$ & $1$\\
$6$ & $144$ & $14$ & $1$\\
$6$ & $280$ & $16$ & $1$\\
$6$ & $628$ & $18$ & $1$\\
$6$ & $512$ & $20$ & $1$\\
$6$ & $432$ & $22$ & $1$\\
$6$ & $900$ & $24$ & $1$\\
\hline
\end{tabular}
\quad
\begin{tabular}{cccc}
\multicolumn{4}{c}{$L_{6,1}$}\\
\hline
$d$ & $n$ & $s$ & $t$\\
\hline
$1$ & $2$ & $1$ & $*$\\
$3$ & $12$ & $4$ & $1$\\
$6$ & $48$ & $6$ & $3$\\
$6$ & $60$ & $8$ & $1$\\
$6$ & $52$ & $7$ & $1$\\
$6$ & $160$ & $12$ & $3$\\
$6$ & $288$ & $14$ & $3$\\
$6$ & $252$ & $16$ & $1$\\
$6$ & $146$ & $13$ & $1$\\
$6$ & $312$ & $20$ & $1$\\
$6$ & $720$ & $22$ & $3$\\
$6$ & $544$ & $24$ & $3$\\
\hline
\end{tabular}
\end{center}

\newpage

\subsection{$d = 7$} We have six $3$-lattices whose shells of norm $3$ are $3$-design, which are $L_{7,3}$, $(\sqrt{3} \mathbb{Z})^7$, $\mathbb{Z}^7$, $\mathbb{Z}^5 \oplus \mathbb{A}_2$, $\mathbb{Z}^3 \oplus (\sqrt{2} \mathbb{Z})^4$, and $L_{7,1}$. A Gram matrix of each lattice is
\begin{equation*}
GM(L_{7,3}) = \left[
\begin{matrix}
 3 & 1 & 1 & 1 & 1 & 1 & 1 \\
 1 & 3 & 1 & 1 & 1 & 1 & 1 \\
 1 & 1 & 3 & 1 & 1 & 1 & 1 \\
 1 & 1 & 1 & 3 & 1 & 1 & 1 \\
 1 & 1 & 1 & 1 & 3 & 1 & -1 \\
 1 & 1 & 1 & 1 & 1 & 3 & -1 \\
 1 & 1 & 1 & 1 & -1 & -1 & 3
\end{matrix}
\right],
 \quad
GM(L_{7,1}) = \left[
\begin{matrix}
 1 & 0 & 0 & 0 & 0 & 0 & 0 \\
 0 & 2 & 1 & 1 & 1 & 1 & 1 \\
 0 & 1 & 2 & 1 & 1 & 1 & 1 \\
 0 & 1 & 1 & 2 & 1 & 1 & 1 \\
 0 & 1 & 1 & 1 & 3 & 0 & 0 \\
 0 & 1 & 1 & 1 & 0 & 3 & 0 \\
 0 & 1 & 1 & 1 & 0 & 0 & 3
\end{matrix}
\right].
\end{equation*}
Here, we have $\det (L_{7,3}) = 64$, $\det ((\sqrt{3} \mathbb{Z})^7) = 2187$, $\det (\mathbb{Z}^7) = 1$, $\det (\mathbb{Z}^5 \oplus \mathbb{A}_2) = 3$, $\det (\mathbb{Z}^3 \oplus (\sqrt{2} \mathbb{Z})^4) = 16$, and $\det (L_{7,1}) = 27$.

Now, the following table is $(d, n, s, t)$-configuration of each shell of norm $m \leqslant 12$ of each lattice:
\begin{center}
\begin{tabular}{r}
\quad\\
\hline
$m$\\
\hline
$1$\\
$2$\\
$3$\\
$4$\\
$5$\\
$6$\\
$7$\\
$8$\\
$9$\\
$10$\\
$11$\\
$12$\\
\hline
\end{tabular}
\quad
\begin{tabular}{cccc}
\multicolumn{4}{c}{$L_{7,3}$}\\
\hline
$d$ & $n$ & $s$ & $t$\\
\hline
\ & \ & \ & \\
\ & \ & \ & \\
$7$ & $56$ & $3$ & $5$\\
$7$ & $126$ & $4$ & $5$\\
\ & \ & \ & \\
\ & \ & \ & \\
$7$ & $576$ & $7$ & $5$\\
$7$ & $756$ & $8$ & $5$\\
\ & \ & \ & \\
\ & \ & \ & \\
$7$ & $1512$ & $11$ & $5$\\
$7$ & $2072$ & $12$ & $5$\\
\hline
\end{tabular}
\quad
\begin{tabular}{cccc}
\multicolumn{4}{c}{$(\sqrt{3} \mathbb{Z})^7$}\\
\hline
$d$ & $n$ & $s$ & $t$\\
\hline
\ & \ & \ & \\
\ & \ & \ & \\
$7$ & $14$ & $2$ & $3$\\
\ & \ & \ & \\
\ & \ & \ & \\
$7$ & $84$ & $4$ & $3$\\
\ & \ & \ & \\
\ & \ & \ & \\
$7$ & $280$ & $6$ & $5$\\
\ & \ & \ & \\
\ & \ & \ & \\
$7$ & $574$ & $8$ & $3$\\
\hline
\end{tabular}
\end{center}

\quad\\

\begin{center}
\begin{tabular}{r}
\quad\\
\hline
$m$\\
\hline
$1$\\
$2$\\
$3$\\
$4$\\
$5$\\
$6$\\
$7$\\
$8$\\
$9$\\
$10$\\
$11$\\
$12$\\
\hline
\end{tabular}
\quad
\begin{tabular}{cccc}
\multicolumn{4}{c}{$\mathbb{Z}^7$}\\
\hline
$d$ & $n$ & $s$ & $t$\\
\hline
$7$ & $14$ & $2$ & $3$\\
$7$ & $84$ & $4$ & $3$\\
$7$ & $280$ & $6$ & $5$\\
$7$ & $574$ & $8$ & $3$\\
$7$ & $840$ & $10$ & $3$\\
$7$ & $1288$ & $12$ & $3$\\
$7$ & $2368$ & $14$ & $3$\\
$7$ & $3444$ & $16$ & $3$\\
$7$ & $3542$ & $18$ & $3$\\
$7$ & $4424$ & $20$ & $3$\\
$7$ & $7560$ & $22$ & $5$\\
$7$ & $9240$ & $24$ & $5$\\
\hline
\end{tabular}
\quad
\begin{tabular}{cccc}
\multicolumn{4}{c}{$\mathbb{Z}^5 \oplus \mathbb{A}_2$}\\
\hline
$d$ & $n$ & $s$ & $t$\\
\hline
$5$ & $10$ & $2$ & $1$\\
$7$ & $46$ & $4$ & $1$\\
$7$ & $140$ & $6$ & $3$\\
$7$ & $330$ & $8$ & $1$\\
$7$ & $592$ & $10$ & $1$\\
$7$ & $786$ & $12$ & $1$\\
$7$ & $1052$ & $14$ & $1$\\
$7$ & $1886$ & $16$ & $1$\\
$7$ & $2710$ & $18$ & $1$\\
$7$ & $2540$ & $20$ & $1$\\
$7$ & $3212$ & $22$ & $1$\\
$7$ & $5740$ & $24$ & $3$\\
\hline
\end{tabular}
\quad
\begin{tabular}{cccc}
\multicolumn{4}{c}{$\mathbb{Z}^3 \oplus (\sqrt{2} \mathbb{Z})^4$}\\
\hline
$d$ & $n$ & $s$ & $t$\\
\hline
$3$ & $6$ & $2$ & $1$\\
$7$ & $20$ & $4$ & $1$\\
$7$ & $56$ & $6$ & $3$\\
$7$ & $126$ & $8$ & $3$\\
$7$ & $232$ & $10$ & $1$\\
$7$ & $392$ & $12$ & $1$\\
$7$ & $576$ & $14$ & $3$\\
$7$ & $756$ & $16$ & $3$\\
$7$ & $1006$ & $18$ & $1$\\
$7$ & $1224$ & $20$ & $1$\\
$7$ & $1512$ & $22$ & $3$\\
$7$ & $2072$ & $24$ & $3$\\
\hline
\end{tabular}
\quad
\begin{tabular}{cccc}
\multicolumn{4}{c}{$L_{7,1}$}\\
\hline
$d$ & $n$ & $s$ & $t$\\
\hline
$1$ & $2$ & $1$ & $*$\\
$3$ & $12$ & $4$ & $1$\\
$7$ & $56$ & $6$ & $3$\\
$7$ & $108$ & $8$ & $1$\\
$7$ & $180$ & $10$ & $1$\\
$7$ & $300$ & $12$ & $1$\\
$7$ & $328$ & $14$ & $1$\\
$7$ & $576$ & $16$ & $1$\\
$7$ & $1018$ & $18$ & $1$\\
$7$ & $760$ & $20$ & $1$\\
$7$ & $984$ & $22$ & $1$\\
$7$ & $2170$ & $24$ & $1$\\
\hline
\end{tabular}
\end{center}

\newpage

\section{Conclusion}
There are some series of $3$-lattices whose shells of norm $3$ are $3$-designs. Before describing them, we want to introduce some known results:

\begin{lemma}[Venkov \cite{V}, Lemma 7.1]
Let $L$ be an even integral lattice of dimension $n \geq 2$ and of minimum $4$, and let $e$ be a minimal vector of $L$. Denote by $p$ the orthogonal projection on the hyperplane $H = e^{\perp}$, put ${L_e}' = \{ x \in L \ | \ (e, x) \equiv 0 \pmod{2} \}$, and let $L_e = p({L_e}')$. Suppose that one of the following two assumptions is verified:
\def\labelenumi{(\arabic{enumi})}
\begin{enumerate}
\item There is $x \in L$ such that $(e, x) \equiv 1 \pmod{2}$;

\item We have $(y, e) \equiv 0 \pmod{2}$ for all $y \in L$, and $L$ contains a vector $x$ such that $(e, x) \equiv 2 \pmod{4}$.
\end{enumerate}
Then, $L_e$ is a odd integral lattice of minimum at least $3$, and we have $\det(L_e) = \det(L)$ under the assumption $(1)$ and $\det(L_e) = \frac{1}{4} \det(L)$ under the assumption $(2)$. We denote $VP_3(L) := L_e$.
\end{lemma}

\begin{theorem}[Pache \cite{P}, parts of Theorem 25 and Proposition 26]\quad
\begin{enumerate}
\item For $n \geqslant 2$, all nonempty shells of $\mathbb{Z}^n$ are spherical $3$-designs.

\item The following shells are spherical $5$-designs:
\begin{align*}
&s_m (\mathbb{Z}^4) & &m = 2 a, \quad a \geqslant 1.\\
&s_m (\mathbb{Z}^7) & &m = 4^a (8 b + 3), \quad a, b \geqslant 0.
\end{align*}

\item For $n \geqslant 2$ and $1 \leqslant m \leqslant 1200$, the nonempty shells of norm $m$ of \ $\mathbb{Z}^n$ are not spherical $5$-designs, except for the above cases.
\end{enumerate}
\end{theorem}

\quad

Now, we have the following three series of $3$-lattices whose shells of norm $3$ are spherical $3$-designs:
\begin{enumerate}
\item $\mathbb{Z}^m$ for $m \geqslant 4$, where only for the case of $m = 7$ we have a spherical $5$-design. (See above theorem)\\

\item $(\sqrt{3} \mathbb{Z})^m$ for $m \geqslant 1$, where there is no spherical $5$-design. (See above theorem)\\

\item $VP_3(\mathbb{A}_2) = \sqrt{3} \mathbb{Z}$, $VP_3(\mathbb{D}_4) = L_{3,3}$, $VP_3(\mathbb{E}_6) = L_{5,3}$, $VP_3(\mathbb{E}_7) = L_{6,3,1}$, and $VP_3(\mathbb{E}_8) = O_7 = L_{7,3}$, where only for the case of $VP_3(\mathbb{E}_8) = O_7$ we have a spherical $5$-design. (cf. Theorem \ref{th-m2-t5})\\
\end{enumerate}

\quad\\

\newpage

\appendix

\section{$3$-lattices of minimum $2$}

\subsection*{d = 2} $\langle 1 \rangle$\\

\noindent \begin{minipage}{1.26in} \noindent $\{ 5, (2,4,3,1), (1,2,1,-) \}$ \par
$\left[

\right]$ \end{minipage}

\quad\\

\subsection*{d = 3} $\langle 4 \rangle$\\

\noindent \begin{minipage}{1.26in} \noindent $\{ 7, (3,6,3,1),$ \par \qquad $(2,6,3,5) \}$ \par
$\left[
%
\right]$ \end{minipage}

\quad\\

\subsection*{d = 4} $\langle 18 \rangle$\\

\noindent \begin{minipage}{1.26in} \noindent $\{ 8, (4,12,5,1),$ \par \qquad $(3,12,4,3) \}$ \par
$\left[
%
\right]$ \end{minipage}

\quad\\

\subsection*{d = 5} $\langle 90 \rangle$

\quad\\
{\bf [No. 1 -- 50]}

\noindent \begin{minipage}{1.26in} \noindent $\{ 8, (5,16,5,1),$ \par \qquad $(4,24,4,5) \}$ \par
{\large $\left[
%
\right]$} \end{minipage}

\quad\\

\subsection*{d = 6} $\langle 591 \rangle$

\quad\\
{\bf [No. 1 -- 100]}

\noindent\begin{minipage}{1.26in} \noindent $\{ 7, (6,32,5,1),$ \par \qquad $(5,40,4,3) \}$ \par
{\Large $\left[
%
\right]$} \end{minipage}

\quad\\
\input{L3-n2-2.tex}
\input{L3-n2-3.tex}

\newpage

\section{$3$-lattices of minimum $3$}

\subsection*{d = 1} $\langle 1 \rangle$\\

\noindent \begin{minipage}{1.57in} \noindent $\{ 3, (1,2,1,-) \}$ \par
$\left[
 3
\right]$ \end{minipage}

\quad\\

\subsection*{d = 2} $\langle 1 \rangle$\\

\noindent \begin{minipage}{1.57in} \noindent $\{ 8, (2,4,3,1) \}$ \par
$\left[
%
\right]$ \end{minipage}

\quad\\

\subsection*{d = 3} $\langle 3 \rangle$\\

\noindent \begin{minipage}{1.57in} \noindent $\{ 16, (3,8,3,{\color{red} \bm{3}}) \}$ \par
$\left[
%
\right]$ \end{minipage}

\quad\\

\subsection*{d = 4} $\langle 9 \rangle$\\

\noindent \begin{minipage}{1.57in} \noindent $\{32, (4,12,3,1) \}$ \par
$\left[
%
\right]$ \end{minipage}

\quad\\

\subsection*{d = 5} $\langle 39 \rangle$

\quad\\
{\bf [No. 1 -- 20]}

\noindent \begin{minipage}{1.57in} \noindent $\{48, (5,20,3,{\color{red} \bm{3}}) \}$ \par
$\left[
%
\right]$ \end{minipage}

\newpage

\subsection*{d = 6} $\langle 247 \rangle$

\quad\\
{\bf [No. 1 -- 100]}

\noindent \begin{minipage}{1.26in} \noindent $\{64, (6,32,3,{\color{red} \bm{3}}) \}$ \par
{\Large $\left[
%
\right]$} \end{minipage}

\quad\\
\subsection*{k = 7} $\langle 2446 \rangle$

\quad\\
{\bf [No. 1 -- 100]}

\noindent \begin{minipage}{1.57in} \noindent $\{ 64, (7,56,3,{\color{red} \bm{5}}) \}$ \par
{\large $\left[
%
\right]$} \end{minipage}

\end{document}